\documentclass[11pt, a4paper, twoside]{amsart}
\usepackage{mathrsfs, tabularx}
\usepackage[all]{xypic}

\usepackage {amsthm, amsmath, amssymb, amscd, float, latexsym, times, epic, eepic}

\theoremstyle{plain}
\newtheorem{theorem}{Theorem}[section]
\newtheorem{lemma}[theorem]{Lemma}
\newtheorem{proposition}[theorem]{Proposition}

\theoremstyle{definition}
\newtheorem{definition}[theorem]{Definition}
\newtheorem{example}[theorem]{Example}

\newcommand {\Set}[1] {\mathbb{#1}}
\newcommand{\setR}[0]{\Set{R}}

\newcommand{\dom}[0]{{\mathscr{D}}}

\newcommand{\semiDistribution}[0]{\Delta}

\newcommand {\proofBox}[0]{\hfill $\Box$ }

\newcommand{\sJ}[0]{j}

\newcommand{\spaceN}[0]{P}

\newcommand{\pd}[2]{\frac{\partial #1}{\partial #2}}

\newcommand{\pdd}[3]{\frac{\partial^2 #1}{\partial #2\,\partial #3}}

\newcommand{\vfield}[1]{{\mathfrak X}( #1)}
\newcommand{\varepsint}[0]{(-\varepsilon, \varepsilon)}

\newcounter{saveenum}

\newcommand{\slaz}[0]{\setminus \{0\}}

\title{A geometric space without conjugate points}

\author[Bucataru]{Ioan Bucataru}
\address{Ioan Bucataru, Faculty of Mathematics, Al.I.Cuza University
B-dul Carol 11, Iasi, 700506, Romania}

\author[Dahl]{Matias F. Dahl}
\address{
Matias F. Dahl, Institute of Mathematics, P.O.Box
1100, 02015 Helsinki University of Technology, Finland }

\date{\today}

\begin{document}
\begin{abstract}
  From a spray space $S$ on a manifold $M$ we construct a new geometric space
  $\spaceN$ of larger dimension with the following properties:
\begin{enumerate}
\item Geodesics in $\spaceN$ are in one-to-one correspondence with parallel
  Jacobi fields of $M$. 
\item $P$ is complete if and only if $S$ is complete.
\item If two geodesics in $P$ meet at one point, the geodesics
  coincide on their common domain, and $P$ has no conjugate points.
\item There exists a submersion 
  that maps geodesics in $P$ into geodesics on $M$.
\end{enumerate}
  Space $\spaceN$ is constructed by first taking two complete lifts
  of spray $S$. This will give a spray $S^{cc}$ on the second iterated
  tangent bundle $TTM$. Then space $\spaceN$ is obtained by restricting
  tangent vectors of geodesics for $S^{cc}$ onto a suitable $(2\dim
  M+2)$-dimensional submanifold of $TTTM$. 
  Due to the last restriction, space $\spaceN$ is not a spray space.
  However, the construction shows that conjugate points can be removed
  if we add dimensions and relax assumptions on the geometric
  structure.
\end{abstract}
\maketitle

\section{Introduction}
\label{sec:intro}
Suppose $S$ is a spray on a manifold $M$. In this paper we show how to
construct a new geometric space $P$ that is based on $S$, but such
that $P$ has no conjugate points.  This is done in three
steps:
\begin{enumerate}
\item We start with a spray $S$ on a manifold $M$.  For example, $S$
  could be the geodesic spray for a Riemannian metric, a Finsler
  metric, or a non-linear connection \cite{BucataruMiron:2007,
    Sakai1992, Shen2001}.
\item Next we take two \emph{complete lifts} of $S$ (see below).  The
  first complete lift gives a spray on $TM$ whose geodesics are Jacobi
  fields on $M$. Similarly, the second complete lift gives a spray
  $S^{cc}$ on $TTM$ whose geodesics can be described as Jacobi fields
  for geodesics for $S^c$. That is, geodesics of $S^c$ describe linear
  deviation of nearby geodesics in $M$, and geodesics of $S^{cc}$
  describe second order deviation of nearby geodesics in $M$.
\item In the last step, we restrict tangent vectors of geodesics of
  $S^{cc}$ onto a submanifold $\semiDistribution\subset TTTM$ that is
  invariant under the geodesic flow of $S^{cc}$. 
  By choosing $\Delta$ in a suitable way, we obtain a space $P$
  where geodesics are in one-to-one correspondence with parallel Jacobi 
  fields in $M$. 
\end{enumerate}
%
In step \emph{(ii)} the original spray $S$ is lifted twice using the
\emph{complete lift}.  Essentially, the complete lift can be seen as a
geometrization of the Jacobi equation. For example, if we start with a
(pseudo-)Riemannian metric $g$ on $M$, the complete lift of $g$ gives
a pseudo-Riemannian metric $g^c$ on $TM$ whose geodesics are Jacobi
fields on $M$.  This means that Jacobi fields on $M$ can be treated as
solutions to a geodesic equation on $TM$, whence there is no need for
a separate Jacobi equation.  In this work we will use the complete
lift of a spray.  For affine sprays, this complete lift was introduced
by A. Lewis \cite{Lewis:2001:GeomMax}. In the Riemannian context, the
complete lift is also known as the \emph{Riemann extension}, and for a
discussion about the complete lift in other contexts, see
\cite{BucataruDahl:2008}.
In step \emph{(ii)}, we need to study sprays on manifolds $M,
TM$, and $TTM$ and also complete lifts of sprays on $M$ and $TM$. To
avoid studying all these cases separately we first study sprays and
complete lifts on iterated tangent bundles of arbitrary order. This is
the topic of Sections
\ref{sec:Notation}-\ref{sec:JacobiFieldsForSpray}.

In step \emph{(iii)} the phase space of spray $S^{cc}$ on $TTM$ is
restricted to a submanifold $\Delta\subset TTTM$. By choosing $\Delta$
suitably, we define a geometry $P$ where geodesics are in
one-to-one correspondence with parallel Jacobi fields (Jacobi fields
of the form $\alpha c'(t) + \beta tc'(t)$ where $c\colon I\to M$ is a
geodesic of $S$). The geometry of sprays that have been restricted in
this way is described in Section \ref{sec:SRTAS}.  Previous work on
sprays with restricted phase space can be found in
\cite{Anastasiei:2001, Lewis:SymmetricProduct, Lewis:Affine1998}.
The space $P$ is constructed and discussed in Section
\ref{sec:SecondNoConjugateSection}. Here we show that $P$ has no
conjugate points. We also show that the canonical submersion
$\pi\colon TTM\to M$ maps geodesics in $P$ into geodesics in
$M$. Hence the geometry of $P$ can be used to study dynamical
properties of $M$.

Let us emphasize that due to the restriction in step \emph{(iii)},
space $\spaceN$ is not a spray space. It seems that to remove
conjugate points, some relaxation of the underlying geometric
structure is needed. For example, in Riemannian geometry the
assumption that a manifold has no conjugate points can have strong
implications.
\begin{enumerate}
\item Suppose $M$ is a $n$-torus with a Riemannian metric.  Then the
no-conjugate assumption implies that $M$ is flat \cite{BuragoIvanov:1994,
Hopf:1948}.
\item Suppose $M$ is a Riemannian manifold such that $M$ is complete,
simply connected, $\dim M\ge 3$, and $M$ is flat outside a compact
set. Then the no-conjugate assumption implies that $M$ is isometric to
$\setR^n$ \cite{Croke:1991}.
\end{enumerate}
See also \cite{Croke:2004, CrokeKleiner:1995, Ruggiero:2007}.  If one
relaxes the assumption on the geometric structure, then the
no-conjugate assumption becomes weaker; on the $2$-torus, there are
non-flat affine connections without conjugate points
\cite{Kikkawa1964}, and on the $n$-torus there are non-flat Finsler
metrics without conjugate points \cite{CrokeKleiner:1995}.

We will not study applications. However, let us note that there are
many problems in both mathematics and physics where a proper
understanding of conjugate points and multi-path phenomena seem to be
important.  For example, in traveltime tomography a typical assumption
is that the manifold has no conjugate points.  See \cite{Croke:2004,
  SteUhl:2005}.
%
Another example is the volume-preserving diffeomorphism group.  This
is a infinite dimensional Riemannian manifold whose geodesics
represent incompressible fluid flows on $M$.  Currently, an active
area of research is studying the exponential map and conjugate points
for this manifold \cite{ArnoldKhesin:1998, Preston:2006}.
As a last example, let us mention geometric optics. Here conjugate
points are problematic since they lead to caustics, where the
amplitude becomes infinite. 

\section{Preliminaries}
\label{sec:Notation}
We assume that $M$ is a smooth manifold without boundary and with
finite dimension $n\ge 1$. By smooth we mean that $M$ is a
topological Hausdorff space with countable base that is locally
homeomorphic to $\setR^n$, and transition maps are $C^\infty$-smooth.
All objects are assumed to be $C^\infty$-smooth on their domains.
 
By $\left(TM, \pi_{0}, M\right)$ we mean the tangent bundle of
$M$. 
For $r\ge 1$, let $T^rM=T\cdots TM$ be the $r$:th \emph{iterated
tangent bundle}, and for $r=0$ let $T^0M=M$. For example, when $r=2$
we obtain the second tangent bundle $TTM$ \cite{Besse:1978,
FisherLaquer:1999}, and in general
   $T^{r+1}M = TT^{r}M$ for $r\ge 0$.

For a tangent bundle $T^{r+1}M$ where $r\ge 0$, we denote
the canonical projection operator by $\pi_{r}\colon T^{r+1}M \to
T^rM$.  Occasionally we also write $\pi_{TTM\to M}$, $\pi_{TM\to M}, 
\ldots$ instead of $\pi_0\circ \pi_1$, $\pi_0, \ldots$.  
Unless otherwise specified, we always use canonical local coordinates 
(induced by local coordinates on $M$) for iterated tangent bundles.
If $x^i$ are local coordinates for $T^rM$ for some $r\ge 0$, we denote
induced local coordinates for $T^{r+1}M$, $T^{r+2}M$, and $T^{r+3}M$
by
\begin{eqnarray*}
& & (x,y), \\
& & (x,y,X,Y), \\
& & (x,y,X,Y,u,v,U,V).
\end{eqnarray*}
As above, we usually leave out indices for local coordinates and write
$(x,y)$ instead of $(x^i, y^i)$. 

For $r\ge 1$, we treat $T^rM$ as a vector bundle over the manifold
$T^{r-1}M$ with the vector space structure induced by projection
$\pi_{r-1}\colon T^rM\to T^{r-1}M$ unless otherwise specified. Thus, if $\{ x^i : i=1, \ldots,
2^{r-1}n \}$ are local coordinates for $T^{r-1}M$, and $(x,y)$ are
local coordinates for $T^rM$, then vector addition and scalar
multiplication are given by
\begin{eqnarray}
\label{eq:TMplus}
  (x,y) + (x, \tilde y) &=& (x,y+\tilde y), \\
\label{eq:mult}
  \lambda \cdot (x,y)  &=& (x,\lambda  y).
\end{eqnarray}
If $x\in T^rM$ and $r\ge 0$ we define
\begin{eqnarray*}
  T_x^{r+1}M &=& \{\xi \in T^{r+1}M : \pi_{r}(\xi)=x\}.
\end{eqnarray*}
For $r\ge 0$, a \emph{vector field} on an open set $U\subset T^rM$
is a smooth map $X\colon U\to T^{r+1}M$ such that $\pi_{r}\circ X =
\operatorname{id}_{U}$.  The set of all vector fields on $U$ is
denoted by $\vfield{U}$.

Suppose that $\gamma$ is a smooth map $\gamma\colon
\left(-\varepsilon, \varepsilon\right)^k \rightarrow T^rM $ where
$k\ge 1$ and $r\ge 0$. Suppose also that $\gamma(t^1,\ldots, t^k) =
(x^i(t^1, \ldots, t^k ))$ in local coordinates for $T^rM$.
Then the \emph{derivative} of $\gamma$ with respect to
variable $t^j$ is the curve $\partial_{t^j}\gamma \colon
\left(-\varepsilon, \varepsilon\right)^k$ $\to T^{r+1}M$ defined by
$\partial_{t^j}\gamma=\left(x^i, {\partial x^i}/{\partial
t^j}\right)$. When $k=1$ we also write $\gamma'=\partial_{t}\gamma$
and say that $\gamma'$ is the \emph{tangent of $\gamma$}.

Unless otherwise specified we always assume that $I$ is an open
interval of $\setR$ that contains $0$, and we do not exclude unbounded
intervals.  If $\phi\colon M\to N$ is a smooth map between
manifolds, we denote the tangent map $TM\to TN$ by $D\phi$, and if
$c\colon I\to M$ is a curve, then
\begin{eqnarray}
\label{curveDiff}
   (\phi\circ c)'(t) = D\phi\circ c'(t), \quad t\in I.
\end{eqnarray}

\subsection{Transformation rules in $T^rM$}
\label{subsec:transformations}
Suppose that $x=(x^i)$ and $\tilde x = (\tilde x^i)$ are overlapping
coordinates for $T^rM$ where $r\ge 0$.  It follows that if $\xi \in
T^{r+1}M$ has local representations $(x, y)$ and $(\tilde x, \tilde
y)$, we have transformation rules
$$
   \tilde x^i = \tilde x^i( x ), \quad
   \tilde y^i = \pd{\tilde x^i}{x^a}(x) y^a.
$$ 
Now $(x, y)$ and $(\tilde x, \tilde y)$ are overlapping coordinates
for $T^{r+1}M$. It follows that if $\xi\in T^{r+2}M$ has local
representations $(x, y,X, Y)$ and $(\tilde x, \tilde y, \tilde X,
\tilde Y)$, we have transformation rules
\begin{eqnarray*}
   \tilde x^i &=& \tilde x^i(x), \\
   \tilde y^i &=& \pd{\tilde x^i}{x^a}(x) y^a, \\
   \tilde X^i &=& \pd{\tilde x^i}{x^a}(x) X^a, \\
   \tilde Y^i &=& \pd{\tilde x^i}{x^a}(x) Y^a + \pdd{\tilde x^i}{x^a}{x^b}(x)  y^a X^b.
\end{eqnarray*}

\section{Lifts on iterated tangent bundles}
\subsection{Canonical involution on $T^rM$}
\label{CanonicalInvolution}
When $r\ge 2$ there are two canonical projections $T^rM\to T^{r-1}M$
given by
\begin{eqnarray}
\label{eq:ProjTrM}
  \pi_{r-1}\colon T^rM\to T^{r-1}M,\quad D\pi_{r-2}\colon T^rM\to T^{r-1}M.
\end{eqnarray}
This means that $T^{r}M$ contains two copies of $T^{r-1}M$, and there
are two ways to treat $T^rM$ as a vector bundle over
$T^{r-1}M$. Unless otherwise specified, we always assume that $T^rM$
is vector bundle $(T^rM, \pi_{r-1}, T^{r-1}M)$, whence the vector
structure of $T^{r}M$ is locally given by equations
\eqref{eq:TMplus}-\eqref{eq:mult}.  However, there is also another
vector bundle structure induced by projection $D\pi_{r-2}\colon
T^rM\to T^{r-1}M$. If $x^i$ are local coordinates for $T^{r-2}M$ and
$(x,y,X,Y)$ are local coordinates for $T^rM$, this structure is given
by
\begin{eqnarray}
\label{eq:vecOther1}
  (x,y,X,Y) + (x, \tilde y, X,\tilde Y) &=& (x,y+\tilde y, X, Y+\tilde Y), \\
\label{eq:vecOther2}
  \lambda \cdot (x,y,X,Y)  &=& (x,\lambda y, X, \lambda Y).
\end{eqnarray}
Next we define the canonical involution $\kappa_r\colon T^rM\to T^rM$
\cite{BucataruDahl:2008}. It is a linear isomorphism between the 
above two vector bundle structures for $T^rM$ defined such that
the following diagram commutes.
\begin{eqnarray*}
\begin{xy}
\xymatrix@C=2.0pc @R=3.0pc{
  T^rM\ar@{<->}[rr]^{\kappa_r} &  & T^rM\ar[dl]^{D\pi_{r-2}} \\
& T^{r-1}M \ar@{<-}[ul]^{\pi_{r-1}}& }
\end{xy}
\end{eqnarray*}
On $TTM$, this involution
map is well known \cite{Besse:1978,
  FisherLaquer:1999,KolarMichorSlovak:1993, Michor1996, Sakai1992}.

\begin{definition}[Canonical involution on $T^rM$]
For $r\ge 2$, the \emph{canonical involution} $\kappa_r \colon T^rM\to
T^rM$ is the unique diffeomorphism that satisfies 
\begin{eqnarray}
\label{KappaDefEq}
  \partial_t \partial_s c(t,s) &=& \kappa_r\circ \partial_s \partial_t c(t,s)
\end{eqnarray}
for all maps $c\colon \varepsint^2 \to T^{r-2}M$. For $r=1$,
we define $\kappa_1=\operatorname{id}_{TM}$.
\end{definition}

Let $r\ge 2$, let $x^i$ be local coordinates for
$T^{r-2}M$, and let $(x,y,X,Y)$ be local coordinates for $T^rM$. Then
\begin{eqnarray*}
  \kappa_r(x,y,X,Y) &=& (x,X,y,Y).
\end{eqnarray*}
For example, in local coordinates for $TTM$ and $TTTM$ we have
\begin{eqnarray*}
  \kappa_2(x,y,X,Y) &=& (x,X,y,Y), \\
  \kappa_3(x,y,X,Y,u,v,U,V) &=& (x,y,u,v,X,Y,U,V).
\end{eqnarray*}
For $r\ge 1$, we have identities
\begin{eqnarray}
\label{eq:noNumberEq}
  \kappa_r^2&=&\operatorname{id}_{T^rM},\\ 
\label{PiKappaXYZ}
  \pi_r\circ D\kappa_r &=& \kappa_r\circ \pi_r, \\ 
\label{commutationRelation}
  D\pi_{r-1} &=&  \pi_{r} \circ \kappa_{r+1}, \\ 
\label{piDpi}
  \pi_{r-1}\circ D\pi_{r-1} &=&  \pi_{r-1}\circ \pi_{r},\\
\label{piDDpi}
  D\pi_{r-1}\circ \pi_{r+1} &=&  \pi_{r}\circ DD\pi_{r-1}, \\
\label{kappaId3}
  DD\pi_{r-1} \circ \kappa_{r+2} &=& \kappa_{r+1} \circ DD\pi_{r-1}, \\
\label{eq:pipikpipi}
  \pi_{r-1}\circ \pi_{r}\circ \kappa_{r+1} &=& \pi_{r-1}\circ \pi_{r}.
\end{eqnarray}

Let us point out that the two
projections in equation \eqref{eq:ProjTrM} are not the only
projections from $T^{r+1}M\to T^rM$. For example, when $r=3$, there
are (at least) $6$ projections $T^3M\to T^2M$; 
$\pi_2$, 
$\kappa_2\circ \pi_2$,
$D\pi_1$, 
$\kappa_2\circ D\pi_1$, 
$DD\pi_0$, and
$\kappa_2 \circ DD\pi_0$.

Let $\gamma_0$ be a curve $\gamma_0\colon I\to T^{r-1}M$ for some $r\ge 1$,
and let
\begin{eqnarray*}
  \vfield{\gamma_0} &=& \{ \eta \colon I\to T^rM : \pi_{r-1}\circ \eta = \gamma_0\}.
\end{eqnarray*}
Elements in $\vfield{\gamma_0}$ are called \emph{vector fields along
$\gamma_0$}, and $\vfield{\gamma_0}$ has a natural vector space
structure induced by the vector bundle structure of $T^rM$ in
equations \eqref{eq:TMplus}-\eqref{eq:mult}.
If $\eta\in \vfield{\gamma_0}$ and $C\in \setR$, then 
\begin{eqnarray}
\label{eq:kappaMult}
  \kappa_{r+1}\circ (C \eta)' &=& C  (\kappa_{r+1}\circ \eta'),
\end{eqnarray}
and if $\eta_1, \eta_2\in \vfield{\gamma_0}$, then
\begin{eqnarray}
\label{eq:kappaPlus}
  \kappa_{r+1}\circ (\eta_1 + \eta_2)' &=& \kappa_{r+1}\circ \eta'_1+ \kappa_{r+1}\circ \eta'_2.
\end{eqnarray}
It follows that $\kappa_{r+1}\circ \partial_t\colon \vfield{\gamma_0}\to 
\vfield{\gamma_0'}$ is a linear map between vector spaces. 

\subsection{Slashed tangent bundles $T^rM\slaz$}
\label{sec:slashTrM}
The \emph{slashed tangent bundle} is the open set in $TM$ 
  defined as
\begin{eqnarray*}
  TM\slaz &=& \{ y\in TM : y\neq 0\}.
\end{eqnarray*}
For an iterated tangent bundle $T^rM$ where $r\ge 2$ we define 
the \emph{slashed tangent bundle} as the open set
\begin{eqnarray*}
  T^rM\slaz &=& \left\{ \xi \in T^rM : (D\pi_{T^{r-1}M\to M})(\xi)\in TM\slaz \right\}.
\end{eqnarray*}
For example, 
\begin{eqnarray*}
   TTM\slaz &=& \{ (x,y,X,Y) \in TTM : X\neq 0\}, \\
   TTTM\slaz &=& \{ (x,y,X,Y,u,v,U,V) \in TTTM : u\neq 0\},
\end{eqnarray*}
where, say, $TTM\slaz = T^2M\slaz$. When $r=0$,  let us also define 
$T^rM\slaz = M$, and for any set $A\subset T^{r}M$ where $r\ge 0$,
let
\begin{eqnarray*}
  A\slaz &=& A \cap T^{r}M\slaz.
\end{eqnarray*}

For $r\ge 1$  we have
\begin{eqnarray}
\label{eq:kappaTrSlaz}
  \kappa_{r+1}(T^{r+1}M\slaz) &=& T(T^rM\slaz), \\
\label{eq:DpiTrMmZ}
  (D\pi_{r-1})(T^{r+1}M\slaz) &=& T^rM \slaz, \\
\label{eq:DkTrMz}
  (D\kappa_{r})(T^{r+1}M\slaz) &=& T^{r+1}M \slaz.
\end{eqnarray}

Before proving these equations, we define the 
\emph{Liouville vector field} 
$E_{r}\in \vfield{T^{r}M}$. For $r\ge 1$, it is given by
\begin{eqnarray*}
  E_r(\xi) &=& \partial_s ( (1+s) \xi)|_{s=0}, \quad \xi\in T^rM.
\end{eqnarray*}
If $r\ge 1$, and $(x,y)$ and $(x,y,X,Y)$ are local
coordinates for $T^{r}M$ and $T^{r+1}M$, respectively, then
\begin{eqnarray*}
  E_r(x,y) &=& (x,y,0,y).
\end{eqnarray*}

Equation \eqref{eq:kappaTrSlaz} follows using equation
\eqref{commutationRelation} and by writing 
\begin{eqnarray}
\label{eq:trivialPiPi}
   \pi_{T^rM\to M} &=& \pi_{T^{r-1}M\to M} \circ \pi_{r-1}, \quad r\ge 1.
\end{eqnarray}
If $r\ge 1$, we have
\begin{eqnarray}
\label{eq:ErId}
  \xi &=& D\pi_{r-1}\circ \kappa_{r+1}\circ E_r(\xi), \quad \xi \in T^rM, 
\end{eqnarray}
and equation \eqref{eq:DpiTrMmZ} follows using equations \eqref{eq:kappaTrSlaz}
and \eqref{eq:trivialPiPi}.
Equation \eqref{eq:DkTrMz} follows using equations \eqref{piDpi} and
\eqref{eq:trivialPiPi}.

\subsection{Lifts for functions}
\label{lifts}
Suppose $f\in C^\infty(M)$ is a smooth function. Then we
can lift $f$ using the \emph{vertical lift} or the \emph{complete
lift} and obtain functions $f^v, f^c\in C^\infty(TM)$ defined by
\begin{eqnarray}
\label{eq:cvl}
  f^v(\xi)=f\circ \pi_0(\xi), \quad f^c(\xi)=df(\xi), \quad \xi\in TM.
\end{eqnarray}
Here $df$ is the exterior derivative of $f$.
In local coordinates $(x,y)$ for $TM$, it follows that
$$
  f^v(x,y)=f(x), \quad f^c(x,y)=\pd{f}{x^i}(x) y^i.
$$
Using these lifts one can define vertical and complete lift
for tensor fields on $M$ of arbitrary order. For a full development
of these issues, see \cite{Yano1973}.

Next we generalize the vertical and complete lifts to functions
defined on iterated tangent bundles $T^rM$ of arbitrary order $r\ge 0$.

\begin{definition}
\label{def:vertLift}
For $r\ge 0$, the 
\emph{vertical lift} of a function $f\in$ $ C^\infty(T^rM\slaz)$ is the
  function $f^v\in C^\infty(T^{r+1}M\slaz)$ defined by
\begin{eqnarray*}
  f^v(\xi) &=&  f\circ \pi_{r} \circ \kappa_{r+1}(\xi), \quad \xi\in T^{r+1}M\slaz.
\end{eqnarray*}
\end{definition}

If $r=0$, Definition \ref{def:vertLift} reduces to lift $f^v$ in
equation \eqref{eq:cvl}, and if $r\ge 1$, equation
\eqref{commutationRelation} implies that $f^v = f\circ D\pi_{r-1}$, 
and equation \eqref{eq:DpiTrMmZ} implies that $f^v$ is smooth.  
For $r\ge 1$, let $x^i$ be local coordinates for
$T^{r-1}M$, and let $(x,y,X,Y)$ be local coordinates for
$T^{r+1}M$. Then
\begin{eqnarray*}
  f^v(x,y,X,Y) &=& f(x,X), \quad f\in C^\infty(T^rM\slaz).
\end{eqnarray*}



\begin{definition}
\label{CompleteLiftFunction}
For $r\ge 0$, the \emph{complete lift} of a function $f\in
C^\infty(T^rM\slaz)$ is the function $f^c\in C^\infty(T^{r+1}M\slaz)$
defined by
\begin{eqnarray*}
  f^c(\xi) &=& (df)\circ \kappa_{r+1}(\xi), \quad \xi\in T^{r+1}M\slaz.
\end{eqnarray*}
\end{definition}

If $r=0$, then Definition \ref{CompleteLiftFunction} reduces to lift
$f^c$ in equation \eqref{eq:cvl}, and if $r\ge 1$, then equation
\eqref{eq:kappaTrSlaz} implies that $f^c$ is smooth.  For $r\ge 1$,
let $x^i$ be local coordinates for $T^{r-1}M$, and let $(x,y,X,Y)$ be
local coordinates for $T^{r+1}M$.  Then
\begin{eqnarray*}
\label{generalFC}
 f^c(x,y,X,Y) &=& \frac{\partial f}{\partial x^a}(x,X)y^a + \frac{\partial f}{\partial y^a}(x,X) Y^a, \quad f\in C^\infty(T^rM\slaz).
\end{eqnarray*}

Taking two complete lifts of $f\in C^\infty(T^rM\slaz)$ yields
\begin{eqnarray}
\label{generalFC2}
  f^{cc}   &=& f^c(x,Y,u,V)  \\
\nonumber
          & & \,\, + \left(\pd{f}{x^a}\right)^c(x,y,u,v) X^a 
            + \left(\pd{f}{y^a}\right)^c(x,y,u,v) U^a,
\end{eqnarray}
where argument $(x,y,X,Y,u,v,U,V) \in T^{r+2}M\slaz$ has been suppressed.

If $f\in C^\infty(T^rM\slaz)$ for some $r\ge 1$, then 
\begin{eqnarray}
\label{eq:vvCom}
  f^{vv} &=& f^{vv}\circ (D\kappa_{r+1}), \\
\label{eq:vcCom}
  f^{vc} &=& f^{cv}\circ (D\kappa_{r+1}), \\
\label{eq:ccCom}
  f^{cc} &=& f^{cc}\circ (D\kappa_{r+1}).
\end{eqnarray}
In Section \ref{sec:JacobiFieldsForSpray} we use these identities to study 
geodesics of iterated complete lifts of a spray. 



\section{Sprays}
A \emph{spray} on $M$ is a vector field $S$ on $TM\slaz$ that
satisfies two conditions. Essentially, these conditions state that
\emph{(i)} integral curves of $S$ are closed under affine
reparametrizations $t\mapsto Ct+t_0$, and \emph{(ii)} an integral
curve of $S$ is of the form $c'\colon I\to TM\slaz$ for a curve
$c\colon I\to M$. Then curve $c\colon I\to M$ is a \emph{geodesic} of
$S$. The motivation for studying sprays is that they provides a
unified framework for studying geodesics for Riemannian metrics,
Finsler metrics, and non-linear connections. See
\cite{BucataruMiron:2007, Sakai1992, Shen2001}.
Next we generalize the definition of a spray to iterated tangent
bundles $T^rM$ for any $r\ge 0$.

\subsection{Sprays on $T^rM$}
\begin{definition}[Spray space] 
\label{def:spray}
Suppose $S$ is a vector field $S\in \vfield{T^{r+1}M\slaz}$ where
  $r\ge 0$.  Then $S$ is a \emph{spray} on $T^{r}M$ if
\begin{enumerate}
\item 
\label{it:SprayI}
$(D\pi_{r})(S)=\operatorname{id}_{T^{r+1}M\slaz}$,
\item 
\label{it:SprayII}
$[E_{r+1}, S] = S$ for Liouville vector field $E_{r+1}\in \vfield{T^{r+1}M}$.
\end{enumerate}
\end{definition}

Let $S$ be a vector field $S\in \vfield{T^{r+1}M\slaz}$ where $r\ge 0$. 
Then condition \ref{it:SprayI} in Definition \ref{def:spray} states that
if $(x,y,X,Y)$ are local coordinates for $T^{r+2}M$, then 
locally 
\begin{eqnarray}
\label{SprayDef}
  S(x,y) &=& \left(x^i,y^i,y^i,-2G^i\left(x,y\right)\right)\\
\nonumber
    &=& \left. y^i \pd{}{x^i}\right|_{(x,y)} - \left. 2 G^i(x,y) \pd{}{y^i}\right|_{(x,y)},
\end{eqnarray} 
where $G^i$ are locally defined functions $G^i\colon T^{r+1}M\slaz \to \setR$.
Condition \ref{it:SprayII} states that functions $G^i$ are 
\emph{positively $2$-homogeneous}; if $(x,y)\in
T^{r+1}M\slaz$, then
\begin{eqnarray*}
  G^i(x,\lambda y) &=& \lambda^2   G^i(x, y), \quad \lambda>0.
\end{eqnarray*}
This is a consequence of Euler's theorem for homogeneous functions
\cite{BCS:2000}.

Conversely, if $S$ is a vector field $S\in \vfield{T^{r+1}M\slaz}$
that locally satisfies these two conditions, then $S$ is a spray on
$T^rM$. Functions $G^i$ in equation \eqref{SprayDef} are called
\emph{spray coefficients} for $S$.

When $r=0$, Definition \ref{def:spray} is equivalent to the usual
definition of a spray \cite{BucataruMiron:2007, Shen2001}. However,
when $r\ge 1$, Definition \ref{def:spray} makes a slightly stronger
assumption on the smoothness of $S$.
Namely, if $r\ge 1$ and $S$ is a spray on $T^rM$ (in the sense of
Definition \ref{def:spray}) then $S$ is smooth on $T^{r+1}M\slaz$, but
if $S$ is a spray on manifold $T^rM$ (in the usual sense) then $S$ is
smooth on $T(T^rM)\slaz$.  Since $T(T^rM)\slaz\supset T^{r+1}M\slaz$,
it follows that if $S$ is a spray on $T^rM$ (in the sense of
Definition \ref{def:spray}), then $S$ is also a spray on manifold
$T^rM$ (in the usual sense).
The stronger assumption on $S$ will be needed in Section
\ref{sec:completeLiftsSprays} to prove that the complete lift of a
spray on $T^rM$ is a spray on $T^{r+1}M$.
In this work we only consider sprays on $T^rM$ that arise from complete
lifts of a spray on $M$. Therefore we do not distinguish between the
weaker and stronger definitions of a spray. These comments
motivate the slightly non-standard terminology in Definition
\ref{def:spray}.

The next proposition shows that a spray on $T^rM$ induces sprays on
all lower order tangent bundles $M, TM, \ldots, T^{r-1}M $.

\begin{proposition} 
\label{prop:projectionOfSpray}
If $S$ 
is a spray on $T^{r+1}M$ where $r\ge 0$, then
\begin{eqnarray*}
   S^\ast &=& (DD\pi_{r}) \circ S \circ \kappa_{r+2}\circ E_{r+1} 
\end{eqnarray*}
is a spray on $T^{r}M$.
\end{proposition}

\begin{proof} 
Equations \eqref{commutationRelation}, \eqref{eq:trivialPiPi}
and equations \eqref{piDDpi}, \eqref{eq:trivialPiPi} 
 imply that maps
\begin{eqnarray*}
  \kappa_{r+2}\circ E_{r+1}\colon T^{r+1}M\slaz &\to& T^{r+2} M \slaz, \\
  DD\pi_r \colon T(T^{r+2}M\slaz) &\to& T(T^{r+1} M \slaz), 
\end{eqnarray*}
are smooth, so $S^\ast \colon T^{r+1}M\slaz\to T(T^{r+1}M\slaz)$ is a
smooth map.
Let $(x,y)$ be local coordinates for $T^{r+1}M$, and let 
$(x,y,X,Y)$ be local coordinates for $T^{r+2}M$. 
Then $S$ can be written as
\begin{eqnarray*}
  S(x,y,X,Y) &=& \left(x,y,X,Y,X,Y,-2G^i(x,y,X,Y), -2H^{i}(x,y,X,Y) \right)
\end{eqnarray*}
for locally defined functions $G^i, H^i\colon T^{r+2}M\slaz \to \setR$
that are positively $2$-homogeneous with respect to $(X,Y)$.
It follows that
\begin{eqnarray*}
  S^\ast (x,y) &=& \left(x,y,y,-2G^i(x,0,y,y)\right),
\end{eqnarray*}
whence $S^\ast$ is a vector field $S^\ast\in \vfield{T^{r+1}M\slaz}$, and
$(D\pi_r)(S^\ast)=\operatorname{id}_{T^{r+1}M\slaz}$. Since functions
$(x,y)\mapsto G^i(x,0,y,y)$ are positively $2$-homogeneous, $S^\ast$
is a spray.
\end{proof}

\subsection{Geodesics on $T^rM$}
Suppose $\gamma$ is a curve $\gamma\colon I \to T^rM$ where $r\ge
0$. Then we say that $\gamma$ is \emph{regular} if $\gamma'(t)\in
T^{r+1}M\slaz$ for all $t\in I$. When $r=0$, this coincides with the
usual definition of a regular curve, and when $r\ge 1$, curve $\gamma$ is
regular if and only if curve $\pi_{T^rM\to M}\circ \gamma \colon I\to M$ is
regular.

\begin{definition}[Geodesic]
\label{def:geodesic}
Suppose $S$ is a spray on $T^rM$ where $r\ge 0$. Then a regular curve
$\gamma \colon I\to T^rM$ is a \emph{geodesic} of $S$ if and only if
\begin{eqnarray*}
  \gamma'' &=& S\circ \gamma'. 
\end{eqnarray*}
\end{definition}

Suppose $S$ is a spray on $T^rM$ and locally $S$ is given by equation
\eqref{SprayDef}. Then a regular curve $\gamma\colon I\to T^{r}M$,
$\gamma=(x^i)$, is a geodesic of $S$ if and only if
\begin{eqnarray}
\label{eq:SGeo}
  \ddot x^i &=& - 2 G^i\circ \gamma'.
\end{eqnarray}

In Definition \ref{def:geodesic} we have defined geodesics on open
intervals.  If $\gamma$ is a curve on a closed interval we say that
$\gamma$ is a geodesic if $\gamma$ can be extended into a geodesic
defined on an open interval.

\section{Complete lifts for a spray}
\label{sec:completeLiftsSprays}
Let $S$ be a spray on $M$. Then the complete lift of $S$ is a spray
$S^c$ on $TM$. That is, if $S$ determines a geometry on $M$, then
$S^c$ determines a geometry on $TM$.  The characteristic feature of
spray $S^c$ is that its geodesics are essentially in one-to-one
correspondence with Jacobi fields of $S$. This correspondence will be
the topic of Section \ref{sec:JacobiFieldsForSpray}. In this section,
we define the complete lift for a spray on an iterated tangent bundle
$T^rM$ of arbitrary order $r\ge 0$. This makes it possible to take
iterated complete lifts; if we start with a spray $S$ on $M$ we can
take iterated lifts $S^c, S^{cc}, S^{ccc}, \ldots$ and lift $S$ onto
an arbitrary iterated tangent bundle.

The definition below for the complete lift of a spray can essentially
be found in \cite[Remark 5.3]{Lewis:2001:GeomMax}. For a further
discussion about related lifts, see \cite{BucataruDahl:2008}.

\begin{definition}[Complete lift of spray] 
\label{def:CompleteLiftOfSpray}
Suppose $S$ is a spray on $T^rM$ for some $r\ge 0$. Then the
\emph{complete lift} of $S$ is the spray $S^c$ on $T^{r+1}M$ defined by
\begin{eqnarray*}
  S^c &=& D\kappa_{r+2} \circ \kappa_{r+3} \circ DS \circ \kappa_{r+2},
\end{eqnarray*}
where $DS$ is the tangent map of $S$, 
\begin{eqnarray*}
  DS\colon T(T^{r+1}M\slaz) &\to& T^2(T^{r+1}M\slaz).
\end{eqnarray*}
\end{definition}

Let us first note that equations \eqref{PiKappaXYZ},
\eqref{commutationRelation}, \eqref{eq:pipikpipi}, and \eqref{eq:trivialPiPi}
imply that
\begin{eqnarray*}
  D\kappa_{r+2}\circ \kappa_{r+3}\colon T^{2}(T^{r+1}M\slaz) &\to & T(T^{r+2}M\slaz)
\end{eqnarray*}
is a smooth map. Thus $S^c$ is a smooth map $T^{r+2}M\slaz\to
T(T^{r+2}M\slaz)$, and by equations \eqref{PiKappaXYZ},
\eqref{commutationRelation}, and \eqref{eq:kappaTrSlaz}, $S^c$ is a
vector field $S^c\in \vfield{T^{r+2}M\slaz}$.  If $S$ is the spray in
equation \eqref{SprayDef}, then locally
\begin{eqnarray}
\label{compSc}
  S^c &=& \left(x,y,X,Y,X,Y,-2(G^i)^v, -2\left(G^i\right)^c \right) \\
\nonumber
      &=& X^i \pd{}{x^i} + Y^i \pd{}{y^i}  -2 (G^i)^v \pd{}{X^i} -2 (G^i)^c
  \pd{}{Y^i},
\end{eqnarray}
and $S^c$ is a spray on $T^{r+1}M$.  

Suppose that $S$ is a spray on $T^rM$ for some $r\ge 0$, and suppose
that $\gamma$ is a regular curve $\gamma\colon I\to T^{r+1}M$,
$\gamma=(x,y)$. Then $\gamma$ is a geodesic of $S^c$ if and only if
\begin{eqnarray}
\label{eq:ScGeo1}
  \ddot x^i &=& - 2 G^i\circ (\pi_r\circ \gamma)', \\
\label{eq:ScGeo2}
   \ddot y^i &=& - 2 (G^i)^c\circ \gamma'.
\end{eqnarray}
It follows that $\pi_{r}\circ \gamma=(x^i)$ is a geodesic of $S$.  In
fact, if $S^\ast$ is the spray in Proposition \ref{prop:projectionOfSpray},
then
\begin{eqnarray}
\label{eq:Srecover}
  S &=&  (S^c)^\ast.
\end{eqnarray}
Thus a spray can always be recovered from its complete lift.
What is more, if $S$ is a spray on $T^{r+1}M$ for $r\ge 0$, then 
$S^{\ast c} = S$ if and only if $S=A^c$ for a spray $A$ on $T^rM$. 

The \emph{geodesic flow} of a spray $S$ is defined as the flow of $S$ as 
a vector field. 

\begin{proposition}[Geodesic flow for the complete lift of a spray]
\label{completeLiftFlow}
Suppose $S$ is a spray on $T^rM$ where $r\ge 0$ and $S^c$ is the
complete lift of $S$. Suppose furthermore that
\begin{eqnarray*}
   \phi\colon \dom(S) \to T^{r+1}M\slaz, \quad    \phi^c\colon \dom(S^c) \to T^{r+2}M\slaz, 
\end{eqnarray*}
are the geodesic flows of sprays $S$ and $S^c$, respectively, with
maximal domains
\begin{eqnarray*}
  \dom(S)\subset T^{r+1}M \slaz \times \setR, \quad \dom(S^c)\subset T^{r+2}M \slaz \times \setR. 
\end{eqnarray*}
Then
\begin{eqnarray}
\label{eq:DSDSC}
\left( (D\pi_r)\times \operatorname{id}_\setR\right) \dom(S^c) &=& \dom(S),
\end{eqnarray}
and
\begin{eqnarray}
\label{eq:SCflow}
  \phi^c_t(\xi) &=& 
  \kappa_{r+2}\circ D\phi_t \circ \kappa_{r+2} (\xi),\quad (\xi,t) \in \dom(S^c),
\end{eqnarray}
where $D\phi_t$ is the tangent map of the map $\xi\mapsto \phi_t(\xi)$
where $t$ is fixed.
\end{proposition}

\begin{proof}
  To prove inclusion ``$\subset$'' in equation \eqref{eq:DSDSC}, let
  $(\xi,t)\in \dom(S^c)$, and let $\gamma\colon I\to T^{r+2}M\slaz$ be an
  integral curve of $S^c$ such that $\gamma(0) = \xi$ and $t\in
  I$. Then
\begin{eqnarray*}
  DD \pi_{r} \circ S^c &=& S \circ D\pi_{r}, 
\end{eqnarray*}
so $D\pi_{r} \circ \gamma\colon I\to T^{r+1}M\slaz$ is an integral curve of
$S$, and $((D\pi_{r})(\xi),t)\in \dom(S)$.
The other inclusion follows similarly since $\gamma'$ is an integral
curve of $S^c$ when $\gamma$ is an integral curve of $S$.

  Suppose $v$ is a curve $v\colon \varepsint \to T^{r+1}M\slaz$, and
  suppose that curve $\xi\colon I\times \varepsint \to T^{r+2}M\slaz$,
\begin{eqnarray}
\label{eq:xitsdef}
  \xi(t,s) &=& \kappa_{r+2}\circ \partial_s\left(\phi_t \circ v(s) \right)
\end{eqnarray}
is defined for some interval $I$ and $\varepsilon>0$. 
For $(t,s)\in I\times \varepsint$ we then have
\begin{eqnarray*}
  S^c\circ \xi(t,s) &=& D\kappa_{r+2} \circ \kappa_{r+3} \circ DS\circ \partial_s\left(\phi_t \circ v(s) \right) \\
             &=& D\kappa_{r+2} \circ \kappa_{r+3} \circ \partial_s (S \circ \phi_t \circ v(s) ) \\
             &=& D\kappa_{r+2} \circ \kappa_{r+3} \circ \partial_s \partial_t (\phi_t \circ v(s) )  \\
             &=& D\kappa_{r+2} \circ \partial_t\partial_s (\phi_t \circ v(s) )  \\
             &=& D\kappa_{r+2} \circ \partial_t \left( \kappa_{r+2}\circ \xi(t,s) \right) \\
             &=& \partial_t \xi(t,s).
\end{eqnarray*}
To prove equation \eqref{eq:SCflow}, let $(\xi_0, t_0)\in \dom(S^c)$.
Let $j(t) = \phi^c_t(\xi_0)$ be the integral curve $j \colon I^\ast
\to T^{r+2}M\slaz$ of $S^c$ with maximal domain $I^\ast\subset \setR$. Then
$t_0\in I^\ast$. 
For a compact subset $K\subset I^\ast$ with $0\in K$ we show that
\begin{eqnarray}
\label{jCompact}
  j(t) = 
  \kappa_{r+2}\circ D\phi_t \circ \kappa_{r+2} (\xi_0),\quad t\in K
\end{eqnarray}
whence equation \eqref{eq:SCflow} follows.
Since $\xi_0\in T^{r+2}M\slaz$, it follows that
$\kappa_{r+2}(\xi_0)=\partial_s v(s)|_{s=0}$ for a curve $v\colon
\varepsint \to T^{r+1}M\slaz$.  Suppose $\tau \in K$. Then $(\xi_0, \tau)\in
\dom(S^c)$, and by equation \eqref{eq:DSDSC}, $(v(0),\tau)\in \dom(S)$.
Since $\dom(S)$ is open 
%
%
\cite{AbrahamMarsden:1994}, there is an open interval $I\ni \tau$ and an
$\varepsilon>0$ such that curve $\xi(t,s)$ in equation
\eqref{eq:xitsdef} is defined on $I\times \varepsint$.
Then, as $K$ is compact, we can shrink $\varepsilon$ and 
assume that $K \subset I$.
Now equation \eqref{jCompact} follows since $\xi(t,0) =
\kappa_{r+2}\circ D\phi_t \circ \kappa_{r+2}(\xi_0)$, $S^c\circ
\xi(t,0) = \partial_t \xi(t,0)$ for $t\in I$, and $\xi(0,0) = \xi_0$.
%
\end{proof}

\section{Jacobi fields for a spray}
\label{sec:JacobiFieldsForSpray}
\begin{definition}[Jacobi field] 
\label{def:JacobiField}
Suppose $S$ is a spray on $T^rM$ where $r\ge 0$, and suppose that
$\gamma \colon I \to T^rM$ is a geodesic of $S$. Then a curve
$J\colon I \to T^{r+1}M$ is a \emph{Jacobi field along} $\gamma$ if
\begin{enumerate}
\item $J$ is a geodesic of $S^c$, 
\item $\pi_r \circ J = \gamma$.
\end{enumerate} 
\end{definition}

In Proposition \ref{CharacterizationOfJacobiFields} we will show that
Definition \ref{def:JacobiField} is equivalent with the usual
characterization of a Jacobi field in terms of geodesic variations. In
view of Proposition \ref{completeLiftFlow}, this should not be
surprising. For example, in Riemannian geometry it is well known that
Jacobi fields are closely related to the tangent map of the
exponential map.
%

\begin{definition}[Geodesic variation] 
  Suppose $S$ is a spray on $T^rM$ where $r\ge 0$, and suppose that
 $\gamma \colon I \to T^rM$ is a geodesic of $S$.  Then a
 \emph{geodesic variation} of $\gamma$ is a smooth map $V\colon
 I\times \varepsint\to T^rM$ such that
\begin{enumerate}
\item $V(t,0)=\gamma(t)$ for all $t\in I$, 
\item $t\mapsto V(t,s)$ is a geodesic for all $s\in \varepsint$.
\end{enumerate}
\end{definition}

Suppose that $I$ is a closed interval. 
Then we say that a curve $J\colon I\to T^rM$ is a Jacobi field if we
can extend $J$ into a Jacobi field defined on an open interval.
Similarly, a map $V\colon I\times \varepsint\to T^{r}M$ is a geodesic
variation if there is a geodesic variation $V^\ast\colon I^\ast \times
(-\varepsilon^\ast, \varepsilon^\ast)\to T^rM$ such that $V=V^\ast$ on
the common domain of $V$ and $V^\ast$ and $I\subset I^\ast$.

The next proposition motivates the above non-standard definition for a
Jacobi field using the complete lift of a spray.

\begin{proposition}[Jacobi fields and geodesic variations]
\label{CharacterizationOfJacobiFields}
Let $S$ be a spray on $T^rM$ where $r\ge 0$, let $J\colon I\to
T^{r+1}M$ be a curve, where $I$ is open or closed, and let
$\gamma\colon I\to T^{r}M$ be the curve $\gamma = \pi_r\circ J$.
\begin{enumerate}
\item If $J$ can be written as
\begin{eqnarray}
\label{eq:JacDef}
  J(t) &=& \left.\partial_s V(t,s)\right|_{s=0}, \quad t\in I
\end{eqnarray}
for a geodesic variation $V\colon I\times \varepsint\to T^rM$, 
then $J$ is a Jacobi field along $\gamma$.
\label{CharacterizationOfJacobiFields:Jac2Sc}
\item If $J$ is a Jacobi field along $\gamma$ and $I$ is compact, then
      there exists a geodesic variation $V\colon I\times \varepsint\to
      T^rM$ such that equation \eqref{eq:JacDef} holds.
\label{CharacterizationOfJacobiFields:Sc2Jac}
\end{enumerate}
\end{proposition}

\begin{proof}
For \ref{CharacterizationOfJacobiFields:Jac2Sc}, let us first
assume that $I$ is open. For $t\in I$ we then have
\begin{eqnarray*}
  S^c\circ \partial_t J(t) &=& D\kappa_{r+2} \circ \kappa_{r+3} \circ DS\circ \kappa_{r+2}\circ \partial_t \partial_s V(t,s) |_{s=0} \\
             &=& D\kappa_{r+2} \circ \kappa_{r+3} \circ \partial_s (S \circ \partial_t V(t,s) )  |_{s=0}\\
             &=& D\kappa_{r+2} \circ \kappa_{r+3} \circ \partial_s \partial_t \partial_t V(t,s) |_{s=0} \\
             &=& D\kappa_{r+2} \circ \partial_t \partial_s \partial_t V(t,s) |_{s=0} \\
             &=& \partial_t \partial_t \partial_s V(t,s) |_{s=0} \\
             &=& J''(t).
\end{eqnarray*}
If $I$ is closed, we can extend $V$ and $J$ so that $I$ is open and 
the result follows from the case when $I$ is open. 

For \ref{CharacterizationOfJacobiFields:Sc2Jac}, we have
 $J'(0)\in T^{r+2}M\slaz$, so we can find a curve $w\colon (-\varepsilon,
 \varepsilon)\to T^{r+1}M\slaz$ such that $\kappa_{r+2} (J'(0)) =
 \partial_s w(s)|_{s=0}$. Then $w(0) = \gamma'(0)$. Since $I$ is
 compact and $\dom(S)$ is open, we can extend $I$ into an open
 interval $I^\ast$ and find an $\varepsilon>0$ such that $V(t,s) =
 \pi_r \circ \phi_t \circ w(s)$ is a map $V\colon I^\ast \times
 (-\varepsilon,\varepsilon)\to T^rM$.
  We have $V(t,0) = \gamma(t)$ for $t\in
 I$, and for each $s\in \varepsint$, the map $t\mapsto V(t,s)$ is a
 geodesic of $S$.  Proposition \ref{completeLiftFlow} and equations
 \eqref{commutationRelation} and \eqref{curveDiff} imply that
 for $t\in I$, 
\begin{eqnarray*}
   J(t) &=& \pi_{r+1} \circ \phi^c_t \circ J'(0) \\
     &=& \pi_{r+1} \circ \kappa_{r+2}\circ D\phi_t \circ \partial_s w(s) |_{s=0} \\
     &=& \partial_s (\pi_r \circ \phi_t \circ w(s) ) |_{s=0}.
\end{eqnarray*}
We have shown that $V$ is a geodesic variation for Jacobi field $J$. 
\end{proof}

Suppose $c\colon I\to M$ is a geodesic for a Riemannian metric, where
$I$ is compact. Then one can characterize Jacobi fields along $c$
using geodesic variations as in Proposition
\ref{CharacterizationOfJacobiFields} \cite{doCarmo:1992}.  Using the
complete lift, we can therefore write the traditional Jacobi equation
in Riemannian geometry as $J'' = S^c\circ J'$.  It is interesting to
note that the derivation of the latter equation only uses
the definition of $S^c$, 
the geodesic equation for $S$, 
the commutation rule \eqref{KappaDefEq} for $\kappa_r$,
and the chain rule in equation \eqref{curveDiff}. 
In particular, there is no need for local coordinates, covariant
derivatives, nor curvature.  For comparison, see the derivations of
the Jacobi equations in Riemannian geometry \cite{Sakai1992}, in Finsler
geometry \cite{BCS:2000}, and in spray spaces \cite{Shen2001}. 
All of these derivations are considerably more involved than the 
proof of Proposition \ref{CharacterizationOfJacobiFields} \emph{(i)}.
For semi-sprays, see also \cite{BucataruMiron:2007} and
\cite{BucataruDahl:2008}.

\subsection{Geodesics of $S^{cc}$}
\label{ex:GeodesicsInScc}
Let $S$ be a spray on $T^rM$ for some $r\ge 0$.  We know that a
regular curve $\gamma\colon I\to T^{r+1}M$, $\gamma=(x,y)$, is a
geodesic of $S^c$ if and only if $\gamma$ locally solves equations
\eqref{eq:ScGeo1}-\eqref{eq:ScGeo2}. Let us next derive corresponding
geodesic equations for spray $S^{cc}$.
 
Let $S$ be given by equation \eqref{SprayDef} in local coordinates
$(x,y)$ for $T^{r+1}M$. Then the complete lift of $S^c$ is the 
spray $S^{cc}$ on $T^{r+2}M$ given by 
\begin{eqnarray*}
   S^{cc} &=& \left(x,y,X, Y, u, v, U, V, u, v, U, V,\right. \\
     & & \,\,\, \left.-2(G^i)^{vv}, -2(G^i)^{cv}, -2(G^i)^{vc}, -2(G^i)^{cc}\right).
\end{eqnarray*}
Suppose $J$ is a regular curve $J\colon I \to T^{r+2}M$,
$J=(x,y,X,Y)$. By equations \eqref{generalFC2} and \eqref{eq:vcCom},
$J$ is a geodesic of $S^{cc}$ if and only if
\begin{eqnarray*}
  \ddot x^i &=&- 2 G^i\circ c', \\
  \ddot y^i &=&- 2 (G^i)^c\circ J_1', \\
  \ddot X^i &=&- 2 (G^i)^c\circ J_2', \\
  \ddot Y^i &=&-2 (G^i)^{cc}\circ J' \\
   &=& - 2 (G^i)^{c}(x^i, Y^i, \dot x^i, \dot Y^i) \\
   & & \,\, - 2\left( 
                  \left(\pd{G^i}{x^a}\right)^c(J'_1) X^a 
                + \left(\pd{G^i}{y^a}\right)^c(J'_1) \dot X^a 
               \right),
\end{eqnarray*}
where curves $c\colon I\to T^rM$, $J_1\colon I\to T^{r+1}M$, and $J_2\colon
I\to T^{r+1}M$ are given by
$$
  c= \pi_{T^{r+2}M\to T^rM}\circ J, \quad
  J_1=\pi_{r+1}\circ J, \quad
  J_2=(D\pi_r)(J),
$$
and in local coordinates $ c= (x^i), J_1= (x^i,y^i)$, and $J_2= (x^i,
X^i).  $

We have shown that if $J$ is a geodesic of spray $S^{cc}$, then $J$
contains two independent Jacobi fields $J_1$ and $J_2$ along $c$.  The
interpretation of this is seen by writing $J=(x,y,X,Y)$ using a
geodesic variation. Then $J_1=(x,y)$ is the base geodesic of $S^c$,
and $J_2=(x,X)$ describes the variation of geodesic $c\colon I\to M$.
A geometric interpretation of components $Y^i$ seems to be more
complicated. For example, $(x,Y)$ does not define a vector field
along $c$. However, for fixed local coordinates, $Y^i$ describe the
variation of the vector components of Jacobi field $J_1=(x,y)$.  If
$J_2=0$, that is, the variation does not vary the base geodesic $c$,
then equations for $Y^i$ simplify and $(x^i, Y^i)$ is a Jacobi field.
In this case, curve $(x^i, Y^i)$ is also independent of local coordinates
(see transformation rules in Section \ref{subsec:transformations}).

\subsection{Iterated complete lifts}
\label{sec:iteratedLiftsForSprays}
Let $S^0$ be a spray on $M$.  For $r\ge 1$, let $S^r$ be the $r$th
iterated complete lift of $S^0$, that is, for $r\ge 1$, let 
\begin{eqnarray*}
  S^{r} &=& (S^{r-1})^c.
\end{eqnarray*}
Then $S^0, S^1, S^2, \ldots$ are sprays on $M, TM, TTM, \ldots$, and
in general, $S^r$ is a spray on $T^rM$.

Equation \eqref{eq:Srecover} shows that each $S^{r}$ contains all
geometry of the original spray $S^0$. A more precise description is
given by equation \eqref{compSc}. It shows that sprays $S^1, S^2,
\ldots$ also contain new geometry obtained from derivatives of spray
coefficients $G^i$ of $S^0$. Namely, the $r$th complete lift $S^{r}$
depends on derivatives of $G^i$ to order $r$. This phenomena can also
be seen from the geodesic flows of higher order lifts. If $\phi$ is
the flow of $S^0$, then up to a permutation of coordinates, the flow
of $S^1$ is $D\phi$, the flow of $S^{2}$ is $DD\phi$, and, in general,
the flow of $S^{r}$ is the $r$th iterated tangent map $D\cdots D\phi$.
This means that the flow of $S^1$ describes the linear deviation of
nearby geodesic of $S$. That is, the flow of $S^1$ describes the
evolution of Jacobi fields.  Similarly, flows of higher order lifts
describe higher order derivatives of geodesic deviations.

\begin{proposition}[New Jacobi fields from old ones]
\label{prop:NewFromOld}
Suppose $S^0, S^1, S^2, \ldots$ are defined as above, and suppose 
that $\sJ \colon I\to T^rM$ is a geodesic for some $S^r$.
\begin{enumerate}
\item \label{prop:NewFromOld:rescale} If $r\ge 0$, $t_0\in \setR$, and
$C>0$, then $\sJ(C t+t_0)$ is a geodesic of $S^r$.
\item \label{prop:NewFromOld:vector} 
If $r\ge 1$ and $k\colon I\to T^rM$ is another geodesic of $S^r$ such 
that $\pi_{r-1}\circ j(t) = \pi_{r-1}\circ k(t)$, then 
$$
  \alpha j + \beta k, \quad \alpha, \beta \in \setR
$$
is a geodesic of $S^r$. 
\item If $r\ge 1$, then $\kappa_r \circ \sJ \colon I \to T^{r}M$ is a
  geodesic of $S^r$.
\label{prop:NewFromOld:III} 
\item If $r\ge 1$, then $\pi_{r-1} \circ \sJ \colon I \to T^{r-1}M$ is
  a geodesic of $S^{r-1}$.
\label{prop:NewFromOld:IV} 
\item If $r\ge 2$, then $(D\pi_{r-2})(\sJ) \colon I \to T^{r-1}M$ is
  a geodesic of $S^{r-1}$.
\label{prop:NewFromOld:phhp} 
\item If $r\ge 0$, then $\sJ' \colon I \to T^{r+1}M$ is a geodesic of
  $S^{r+1}$. \label{prop:NewFromOld:I} 
\item If $r\ge 0$, then $t\sJ'(t) \colon I \to T^{r+1}M$
  is a geodesic of $S^{r+1}$. \label{prop:NewFromOld:Ib} 
\item If $r\ge 1$, then $E_r\circ \sJ \colon I \to T^{r+1}M$ is a
  geodesic of $S^{r+1}$. \label{prop:NewFromOld:IIX} 
\end{enumerate}
\end{proposition}

\begin{proof} 
  Properties \ref{prop:NewFromOld:rescale},
  \ref{prop:NewFromOld:vector}, and \ref{prop:NewFromOld:IV} follow
  using equations \eqref{eq:SGeo}, \eqref{eq:ScGeo1}, and \eqref{eq:ScGeo2}.
  Properties \ref{prop:NewFromOld:I}, \ref{prop:NewFromOld:Ib}, 
  and \ref{prop:NewFromOld:IIX} follow 
  by locally studying geodesic variations
\begin{eqnarray*}
  V(t,s) &=& j(t+s), \\
  V(t,s) &=& j((1+s)t), \\
  V(t,s) &=& (1+s) j(t),
\end{eqnarray*}
and using Proposition \ref{CharacterizationOfJacobiFields} 
\emph{(i)}.
Property \ref{prop:NewFromOld:III} follows using geodesic equations
for $S^{cc}$ in Section \ref{ex:GeodesicsInScc} and equation
\eqref{eq:ccCom}.  Property \ref{prop:NewFromOld:phhp} follows using
equation \eqref{commutationRelation}.
\end{proof}

\subsection{Conjugate points}
\label{sec:ConjugatePointsForSc}
Suppose $S$ is a spray on $T^rM$ for some $r\ge 0$.  If $a,b$ are
distinct points in $T^rM$ that can be connected by a geodesic $\gamma
\colon [0,L] \to T^rM$, then $a$ and $b$ are \emph{conjugate points} if
there is a Jacobi field $J\colon[0,L]\to T^{r+1}M$ along $\gamma$ that
vanishes at $a$ and $b$, but $J$ is not identically zero (with respect
to vector space structure in equations
\eqref{eq:TMplus}-\eqref{eq:mult}).

The next proposition shows that $S$ has conjugate points if and only
if $S^c$ has conjugate points. Thus the complete lift alone does not
remove conjugate points. 

\begin{proposition}[Conjugate points and complete lift] 
\label{prop:conjuagePoints}
Suppose $S$ is a spray on $T^rM$ for some $r\ge 0$. 
\begin{enumerate}
\item 
\label{it:kajsdsa}
If $a,b \in T^rM$ are conjugate points for $S$, then zero vectors in
$T_a^{r+1}M$ and $T_b^{r+1}M$ are conjugate points for $S^c$.
\item 
\label{it:kajsdsa1a}
If $a,b \in T^rM$ are conjugate points for $S$, then there are
non-zero conjugate points in $T^{r+1}_aM$ and $T^{r+1}_bM$ for $S^c$.
\item 
\label{it:kajsdsa2}
If $a,b\in T^{r+1}M$ are conjugate points for $S^c$, then $\pi_r(a),
\pi_r(b)$ are conjugate points for $S$.
\end{enumerate}
\end{proposition}

\begin{proof} 
%
%
  For property \ref{it:kajsdsa}, suppose $J\colon [0,L]\to T^{r+1}M$
  is a Jacobi field of $S$ that shows that $a$ and $b$ are conjugate
  points.  Then the claim follows by studying Jacobi field
  $E_{r+1}\circ J$.
%
%
  For property \ref{it:kajsdsa1a}, suppose that $J\colon [0,L]\to
  T^{r+1}M$ is as in \ref{it:kajsdsa}, and let $\gamma\colon [0,L]\to
  T^rM$ be the geodesic $\gamma=\pi_r\circ J$ for $S$. We will show
  that $\gamma'(0),\gamma'(L)\in T^{r+1}M\slaz$ are conjugate points
  for $S^c$.  
This follows by considering Jacobi field $j\colon [0,L]\to
  T^{r+2}M$,
\begin{eqnarray*}
  \sJ(t) &=& \partial_s \left( \gamma'(t) + sJ(t)\right) |_{s=0}.
\end{eqnarray*}
The claim follows since $\sJ$ vanishes
at $0$ and $L$, but $\sJ$ is not identically zero.
%
%
For property \ref{it:kajsdsa2}, suppose $J\colon [0,L]\to T^{r+2}M$ is
a Jacobi field of $S^c$ that shows that $a$ and $b$ are conjugate
points.  Then $J$ is a geodesic of $S^{cc}$, and locally $J$ satisfy
equations in Section \ref{ex:GeodesicsInScc}.  If Jacobi field $j =
(D\pi_{r})(J)$ does not vanish identically, the claim follows.
Otherwise $(D\pi_{r})(J)$ vanishes identically, and the result
follows by the last comment in Section \ref{ex:GeodesicsInScc}.
\end{proof}

\section{Sprays restricted to a semi-distribution}
\label{sec:SRTAS}
From a spray $S$ on $M$ one can construct a new geometric space by
restricting the spray to a geodesically invariant distribution
$\Delta\subset TM$. This is done by requiring that all geodesics are
tangent to the distribution.  
For example, geodesics in Euclidean space $\setR^3$ can in this way be
constrained to $xy$-planes.
See
\cite{Anastasiei:2001, Lewis:SymmetricProduct, Lewis:Affine1998}.

In this section we study a slightly more general geometry, where one
can not only restrict possible directions, but also basepoints for
geodesics.  For example, geodesics in $\setR^3$ can in this way be
constrained to one line or one plane.  For a spray on $T^rM$, this is
done by requiring that geodesics are tangent to a suitable
geodesically invariant submanifold $\Delta\subset T^{r+1}M$. Such a
submanifold will be called a \emph{semi-distribution} and the
restricted geometry will be called a \emph{sub-spray}.  There does not
seem to be any work on this type of geometry.  The terms
semi-distribution and sub-spray neither seem to have been used before.

\begin{definition}
\label{def:SemiDistribution}
A set $\semiDistribution\subset T^{r+1}M$ where $r\ge 0$ is a
\emph{semi-distribution on $T^rM$} if 
\begin{enumerate}
\item $\pi_r(\Delta)$ is a submanifold in $T^rM$.
\item $B=\pi_r\circ \kappa_{r+1}(\Delta)$ is a submanifold in $T^rM$.
\item There is a $k\ge 1$ such that every $b \in B$ has an
  open neighbourhood $U\subset B$, and there are $k$ maps
  $V_1, \ldots, V_k\colon U\to T^{r+1}M$ such that
\begin{enumerate}
\item[\emph{(a)}] $\pi_r\circ V_i= \iota$ for $i=1,\ldots, k$, where $\iota$ is inclusion $U\hookrightarrow T^rM$,
\item[\emph{(b)}] $V_i$ are pointwise linearly
  independent, 
\item[\emph{(c)}] for all $u \in U$ we have
\begin{eqnarray*}
  \kappa_{r+1}(\Delta)\cap \pi_r^{-1}(u) &=& \operatorname{span}\{ V_1(u), \ldots, V_k(u)\}.
\end{eqnarray*}
\end{enumerate}
(In \emph{(b)} and \emph{(c)}, the linear structure of $T^{r+1}M$
is with respect to equations \eqref{eq:TMplus}-\eqref{eq:mult}.)
\end{enumerate}
We say that $k$ is the \emph{rank} of $\Delta$ and write $\operatorname{rank} \Delta=k$.
\end{definition}
In condition \emph{(ii)}, $B=\pi_0(\Delta)$ when $r=0$, and $B=
(D\pi_{r-1})(\Delta)$ when $r\ge 1$. Thus, if $r=0$ and
$\pi_0(\Delta)=M$, a semi-distribution is a distribution in the usual
sense.
Condition \emph{(iii)} states that there is a $k$ dimensional vector
space associated to each $b\in B$, and $1\le k \le 2^r \dim M$.  When
$r= 0$, the structure of these vector spaces in $\Delta$ is given by
equations \eqref{eq:TMplus}-\eqref{eq:mult}, and when $r\ge 1$, the
structure is given by equations
\eqref{eq:vecOther1}-\eqref{eq:vecOther2}.
The next example motivates the use of vector space structure in
equations \eqref{eq:vecOther1}-\eqref{eq:vecOther2} when $r\ge
1$. Namely, these equations describe the natural vector space
structure for tangents to Jacobi fields.

\begin{example}
\label{ex:ScVStr}
Let $S$ be a spray on $T^rM$ for some $r\ge 0$, let $\gamma\colon I\to
T^rM$ be a geodesic of $S$, and let $\vfield{\gamma}$ be the set of
vector fields along $\gamma$ with the vector space structure defined
by equations \eqref{eq:TMplus}-\eqref{eq:mult}.
Furthermore, let $J_1, J_2\in \vfield{\gamma}$ be Jacobi fields along
$\gamma$, such that locally $\gamma = (x), J_1=(x,y),$ and
$J_2=(x,z)$. For $\alpha, \beta\in \setR$ we then have
\begin{eqnarray*}
   \alpha J_1 + \beta J_2 &=& (x,\alpha y  + \beta z), \\
   (\alpha J_1 + \beta J_2)' &=& (x,\alpha y  + \beta z, \dot x,\alpha \dot y  + \beta \dot z) \\
   &=& \alpha \cdot J_1'  + \beta  \cdot J_2',
\end{eqnarray*}
where on the last line, $+$ and $\cdot$ are as in equations
\eqref{eq:vecOther1}-\eqref{eq:vecOther2}. 
Thus, if we define the vector space structure for Jacobi fields
by equations \eqref{eq:TMplus}-\eqref{eq:mult}, then the natural
vector structure for tangents (and initial values) 
is given by
equations \eqref{eq:vecOther1}-\eqref{eq:vecOther2}.  
%
%
On the other hand, the multiplication operator in equation
\eqref{eq:mult} appears naturally when reparametrizing a curve.  If
$J\colon I\to T^rM$ is a curve for $r\ge 0$, and $j(t)=J(Ct + t_0)$,
then $j'(t)=C\cdot J'(Ct+t_0)$, where $\cdot$ is as in equation
\eqref{eq:mult}.  \proofBox
\end{example}

\begin{proposition}
\label{prop:SemiDistributionChar}
Suppose $\Delta$ is a
semi-distribution on $T^rM$ and $B=\pi_r\circ
\kappa_{r+1}(\Delta)$. Then $\Delta$ is a sub-manifold in $T^{r+1}M$
and
\begin{eqnarray*}
   \dim \Delta &=& \dim B + \operatorname{rank} \Delta.
\end{eqnarray*}
\end{proposition}

The proof of Proposition \ref{prop:SemiDistributionChar} follows by
setting $A=\kappa_{r+1}(\Delta)$ in the lemma below. We also use this
lemma to prove Proposition \ref{prop:step3a}.

\begin{lemma}
\label{lemma:subman}
Suppose $A$ is a subset $A\subset T^{r+1}M$ for some $r\ge 0$ such that
\begin{enumerate}
\item $\pi_r(A)$ is a submanifold in $T^rM$.
\item There is a $k\ge 1$ such that every $b \in \pi_r(A)$ has an open
  neighbourhood $U\subset \pi_r(A)$, and there are $k$ maps $V_1, \ldots,
  V_k\colon U\to T^{r+1}M$ such that
\begin{enumerate}
\item[\emph{(a)}] $\pi_r\circ V_i= \iota$ for $i=1,\ldots, k$,  where $\iota$ is inclusion $U\hookrightarrow T^rM$,
\item[\emph{(b)}] 
$V_1, \ldots, V_k$ are pointwise linearly independent in $U$, 
\item[\emph{(c)}] for all $u \in U$ we have
\begin{eqnarray*}
  A \cap \pi^{-1}_r(u) &=& \operatorname{span}\{ V_1(u), \ldots, V_k(u)\}.
\end{eqnarray*}
\end{enumerate}
(In (b) and (c), the linear structure of $T^{r+1}M$ is 
with respect to equations \eqref{eq:TMplus}-\eqref{eq:mult}.)
\end{enumerate}
Then $A$ is a submanifold of $T^{r+1}M$ of dimension $\dim \pi_r(A) + k$.
Moreover, if we can assume that $U=\pi_r(A)$, then $A$ is diffeomorphic
to $\pi_r(A)\times \setR^k$. 
\end{lemma}

\begin{proof} Let $\xi \in A$. Then $\pi_r(\xi)$ has an open
  neighbourhood $U\subset \pi_r(A)$ with $k$ maps $V_1, \ldots,
  V_k\colon U\to T^{r+1}M$ such that \emph{(a)}, \emph{(b)}, and
  \emph{(c)} hold.  By possibly shrinking $U$ we can find maps
  $V_{k+1}, \ldots, V_N\colon U\to T^{r+1}M$, where $N=\dim
  T_{\pi_r(\xi)}^{r+1}M$, such that 
$\pi_r\circ V_i= \iota$ for $i=1,\ldots, N$, and
for all $u \in U$, 
\begin{eqnarray*}
  \pi^{-1}_r(u) &=& \operatorname{span}\{ V_1(u), \ldots, V_N(u)\}.
\end{eqnarray*}
%
Let $f$ be the diffeomorphism $f\colon U\times \setR^N\to \pi_r^{-1}(U)$ 
defined as
\begin{eqnarray*}
   f(u,\alpha_1, \ldots, \alpha_N) &=& 
     \alpha_1 V_1(u) + \cdots + \alpha_N V_N(u).
\end{eqnarray*}
Let $g\colon U\times \setR^k\to \pi_r^{-1}(U)$ be the restriction of
$f$ onto $U\times \setR^k$. Then $g$ is a smooth injection and
immersion such that $g(U\times \setR^k)=A\cap \pi_r^{-1}(U)$, and map
$f^{-1}\circ g\colon U\times \setR^k\to U\times \setR^N$ is the
inclusion $(u,\alpha_1, \ldots, \alpha_k)\mapsto (u,\alpha_1, \ldots,
\alpha_k,0,\ldots,0)$. 
Since a closed set in a compact Hausdorff space is compact,
$f^{-1}\circ g$ is proper.  Thus $g$ is proper, and the claim follows
from the following result: If $f\colon M\to N$ is a smooth map between
manifolds that is proper, injective, and an immersion, then $f(M)$ is
a submanifold in $N$ of dimension $\dim M$, and $f$ restricts to a
diffeomorphism $f\colon M\to F(M)$. See results 7.4, 8.3, and 8.25 in
\cite{Lee:abc}.
\end{proof}

\subsection{Geodesics in a sub-spray}
\label{sec:GeoSubX}
\begin{definition}[Geodesically invariant set] 
  Let $S$ be a spray on $T^rM$ where $r\ge 0$. Then a set
  $\semiDistribution\subset T^{r+1}M$ is a \emph{geodesically
    invariant set for $S$} provided that:
\begin{enumerate}
\item[] If $\gamma\colon I\to T^rM$ is a geodesic of $S$ with
  $\gamma'(t_0)\in \semiDistribution$ for some $t_0\in I$, then
  $\gamma'(t)\in \semiDistribution$ for all $t\in I$.
\end{enumerate}
\end{definition}

\begin{definition}[Sub-spray] 
\label{def:geoSubSpray}
Suppose $S$ is a spray on $T^rM$ for some $r\ge 0$, and
$\semiDistribution$ is a geodesically invariant semi-distribution on
$T^rM$.
Then we say that triple $\Sigma = (S, T^rM, \semiDistribution)$ is a
\emph{sub-spray}. A curve $\gamma\colon I\to T^rM$ is a \emph{geodesic
  in $\Sigma$} if
\begin{enumerate}
\item $\gamma\colon I\to T^rM$ is a geodesic of $S$,
\item $\gamma'(t_0)\in \semiDistribution$ for some $t_0\in I$ (whence $\gamma'(t)\in \Delta$ for all $t\in I$).
\label{subSprayCond2}
\end{enumerate}
\end{definition}

By taking $\semiDistribution = T^{r+1}M$, we may treat any spray as a
sub-spray. 
Let us also note that if $\Delta\subset T^{r+1}M\slaz$ where
$r\ge 0$, then
$$
  \pi_r(\Delta)\subset T^rM, \quad \pi_r\circ \kappa_{r+1} (\Delta) \subset T^rM\slaz.
$$

Let $\Sigma = (S, T^rM, \semiDistribution)$ be a sub-spray for some
$r\ge 0$.  
Then
\begin{eqnarray*}
  \semiDistribution\slaz &=& \left\{ \gamma'(0) 
            : \, \gamma\colon
 \varepsint\to T^rM \mbox{ is a geodesic in $\Sigma$}\,\right\}, \\
  \pi_{r}(\semiDistribution\slaz) &=& \left\{ \gamma(0) 
            : \, \gamma\colon \varepsint\to T^rM \mbox{ is a geodesic in $\Sigma$}\,\right\}.
\end{eqnarray*}
In other words, a vector $\xi\in T^{r+1}M$ is in
$\semiDistribution\slaz$ if and only if there is a geodesic in
$\Sigma$ whose tangent passes through $\xi$, and a point $x\in T^rM$
is in $\pi_{r}(\semiDistribution\slaz)$ if and only if there is a geodesic
in $\Sigma$ that passes through $x$.  We therefore say that
$\Delta\slaz$ is \emph{phase space} for $\Sigma$, and
$\pi_r(\Delta\slaz)$ is \emph{configuration space} for $\Sigma$.
When $r\ge 1$, the set $B=(D\pi_{r-1})(\Delta)$ satisfies
\begin{eqnarray*}
  B\slaz &=& \left\{ (\pi_{r-1}\circ \gamma)'(0) 
            : \, \gamma\colon
 \varepsint\to T^rM \mbox{ is a geodesic in $\Sigma$}\,\right\}.
\end{eqnarray*}
and we can interpret $B\slaz$ as phase space of geodesics in $\Sigma$
that have been projected onto $T^{r-1}M$. 

\begin{example}[Geodesics through a point]
  Let $\Sigma =(S, T^rM, \Delta)$ be a sub-spray for some $r\ge 0$,
  and let $z\in \pi_r(\Delta\slaz)$ be a point in configuration space.  Then
  the set
\begin{eqnarray*}
  \Delta(z) &=& \Delta\cap (T_z^{r+1}M \slaz)
\end{eqnarray*}
parametrizes initial values for geodesics that pass through
$z$. Let us study the structure and the degrees of freedom for
$\Delta(z)$.

When $r=0$, the structure of $\Delta(z)$ is easy to understand; the
set $\Delta(z)$ is a punctured vector subspace of $T_zM$ whose
dimension is the rank of $\Delta$.

When $r\ge 1$, the structure of $\Delta(z)$ becomes more complicated.
For example, in Section \ref{sec:SecondNoConjugateSection}, we
construct a sub-spray where configuration space 
and phase space 
are diffeomorphic, and $\Delta(z)$ contains only one vector.  To
understand this, let us assume that $\Delta$ is represented in
canonical local coordinates $(x,y,X,Y)$ for $T^{r+1}M$. That is, we
here only consider coordinates $(x,y,X,Y)$ that belong to
$\Delta$. Then coordinates $(x,y,X,Y)$ have $\dim \Delta = \dim B +
\operatorname{rank} \Delta$ degrees of freedom. Coordinates $(x,X)$
represent submanifold $B$. They have $\dim B$ degrees of freedom, and
once $(x,X)\in B$ is fixed, coordinates $(y,Y)$ parametrize the
$\operatorname{rank}\Delta$ dimensional vector space associated with
$(x,X)$.  If $z=(x_0,y_0)$, then geodesics that pass through $z$ are
parametrized by $(x_0, y_0, X, Y)$, but very little can be said about
possible values for $(X,Y)$.  Coordinates $(x,X)$ have $\dim B$
degrees of freedom, but we do not know how these divide between $x$-
and $X$-coordinates.  Similarly, coordinates $(y,Y)$ have
$\operatorname{rank} \Delta$ degrees of freedom, but we do not know
how these divide between $y$- and $Y$-coordinates.  \proofBox
\end{example}

The next proposition shows that geodesics in a sub-spray on $T^rM$
have a linear structure when $r\ge 1$, but geodesics are not necessarily 
invariant under affine reparametrizations.

\begin{proposition}
\label{prop:NewFromOld:subspray}
Let $\Sigma=(S,T^rM, \Delta)$ be a sub-spray where $r\ge 0$. 
\begin{enumerate}
\item Suppose that $j\colon I\to T^rM$ is a geodesic in $\Sigma$.
  If $t_0\in \setR$, and $C>0$, then $j(C t+t_0)$ is a
  geodesic in $\Sigma$ if $r=0$ or $C=1$. 
\item Suppose that $r\ge 1$. If $j,k\colon I\to T^rM$ are geodesics
  in $\Sigma$ such that $\pi_{r-1}\circ j(t) = \pi_{r-1}\circ k(t)$,
  then
$$
  \alpha j + \beta k, \quad \alpha,\beta \in \setR
$$
is a geodesic in $\Sigma$. 
\end{enumerate}
\end{proposition}

\begin{proof} 
Property \emph{(i)} follows since reparametrizations scale tangent
vectors as in equation \eqref{eq:mult}, and this 
multiplication is only compatible with the vector structure of $\Delta$
when $r=0$ or $C=1$. 
Property \emph{(ii)} follows using 
equations \eqref{eq:kappaMult}-\eqref{eq:kappaPlus}.
\end{proof}

\subsection{Jacobi fields for a sub-spray}
Proposition \ref{CharacterizationOfJacobiFields} shows that for
sprays, Jacobi fields on compact intervals can be characterized using
geodesic variations.  For sub-sprays, we take this characterization as
the definition of a Jacobi field. 

\begin{definition}[Jacobi field in a sub-spray] 
\label{def:JacobiFieldSubSpray:SS}
Let $\gamma\colon I\to T^rM$ be a geodesic in a sub-spray $\Sigma =
(S, T^rM, \semiDistribution)$ where $r\ge 0$. Suppose that $J \colon K
\to T^{r+1}M$ is a curve where $K\subset I$ is compact, and $V$ is a
map $V\colon I\times \varepsint\to T^rM$ such that
\begin{enumerate}
\item $t\mapsto V(t,s)$, $t\in I$ is a geodesic in sub-spray $\Sigma$
  for all $s\in \varepsint$,
\item $V(t,0) = \gamma(t)$ for $t\in I$, 
\item $J(t) = \partial_s V(t,s)|_{s=0}$ for $t\in K$.
\end{enumerate}
Then $J\colon K\to T^{r+1}M$ is a \emph{Jacobi field along $\gamma$}.
\end{definition}

By Proposition \ref{CharacterizationOfJacobiFields} \emph(ii), a
Jacobi field for a sub-spray $(S, T^rM, \Delta)$ is a Jacobi field for
the spray $S$.  The converse also holds when $\Delta = T^{r+1}M$.

\section{A sub-spray for parallel Jacobi fields}
\label{sec:SecondNoConjugateSection}
This section contains the main results of this paper. We construct a
sub-spray $\spaceN$ whose geodesics are in one-to-one correspondence
with parallel Jacobi fields, and study its geodesics. 


\begin{definition}[Parallel Jacobi field]
\label{ex:parJac}
Let $S$ be a spray on $M$, and let $J\colon I\to TM$ be a curve. Then
$J$ is called a \emph{parallel Jacobi field for $S$} if there are
$\alpha,\beta\in \setR$, and a geodesic $c\colon I\to M$
such that
\begin{eqnarray}
\label{eq:ParJacFie}
  J(t) &=& \alpha c'(t) + \beta tc'(t), \quad t\in I.
\end{eqnarray}
\end{definition}

If $I$ is closed we say that a curve $J\colon I\to TM$ is a parallel
Jacobi field if $J$ can be extended into a parallel Jacobi defined on
an open interval. Proposition \ref{prop:NewFromOld} shows that a
parallel Jacobi field is a Jacobi field.

\begin{lemma}
\label{lemma:extensionResult}
Suppose $J \colon I\to TM$ is a parallel Jacobi field.
\begin{enumerate}
\item If $C>0$ and $t_0\in \setR$, then $J(Ct + t_0)$ is a parallel
  Jacobi field.
\item $J$ can be extended to the maximal domain of geodesic
  $c=\pi_{TM\to M}\circ J$, and the extension is a parallel Jacobi
  field.
\end{enumerate}
\end{lemma}




To construct sub-spray $P$, let $S$ be a spray on a manifold $M$, let
$S^{cc}$ be the second complete lift of $S$, and let
$\semiDistribution$ be the geodesically invariant semi-distribution on
$TTM$ defined in Proposition \ref{prop:step3a}. Then we define
sub-spray $P$ as
\begin{eqnarray*}
  \spaceN &=& (S^{cc},TTM, \semiDistribution).
\end{eqnarray*}

\begin{proposition} 
\label{prop:step3a}
Suppose $S$ is a spray on a manifold $M$, and let $\semiDistribution$ be the set
\begin{eqnarray*}
  \semiDistribution &=& \left\{ (\kappa_2\circ J')'(0): J\colon \varepsint \to TM \mbox{\, is parallel Jacobi field for $S$} \right\}.
\end{eqnarray*}
Then 
\begin{enumerate}
\item $\semiDistribution\subset TTTM\slaz$, 
\item $\semiDistribution$ is a geodesically invariant
  semi-distribution on $TTM$ of rank $2$,
\item phase space $\Delta$ and configuration space $\pi_2(\Delta)$ are
  diffeomorphic.
\end{enumerate}
\end{proposition}

\begin{proof}
  Let us first note that $\semiDistribution$ consists
  of points 
\begin{eqnarray*}
  \left(x(0), 
   \dot x(0), 
   \alpha \dot x(0), 
   \alpha \ddot x(0) + \beta \dot x(0), \right.
   \quad\quad\quad\quad\quad\quad
   \\
\left.   
   \dot x(0), 
   \ddot x(0), 
   \alpha \ddot x(0) + \beta \dot x(0), 
   \alpha \dddot x(0) +  2\beta \ddot x(0)  
  \right), 
\end{eqnarray*}
where $\alpha, \beta\in \setR$ and $c\colon \varepsint\to M$ is a geodesic
$c(t)=(x(t))$. 
By Lemma \ref{lemma:subman} (and by the result used to 
prove Lemma \ref{lemma:subman}), it follows that sets
\begin{eqnarray*}
 \pi_2(\Delta) &=& \{\alpha S(y) + \beta E_1(y) : y\in TM\slaz, \alpha, \beta\in \setR \}, \\
  (D\pi_1)(\Delta) &=& \{ S(y) : y\in TM\slaz\}, 
\end{eqnarray*}
are submanifolds in $TTM$ diffeomorphic to $TM\slaz\times \setR^2$
and $TM\slaz$, respectively.
Let $\iota$ be the inclusion $B\hookrightarrow TTM$,
where $B=(D\pi_1)(\Delta)$, and let $\widehat S$ be the 
diffeomorphism $\widehat S\colon TM\slaz \to B$ such that $S= \iota \circ \widehat S$ and $\widehat S^{-1} = \pi_1\circ \iota$. 
By the geodesic equation for $S^c$ and the definition of $S^c$ 
it follows that
\begin{eqnarray*}
 \kappa_3(\Delta) 
    &=& (DS) (\pi_2(\Delta)) \\
    &=& \{\alpha V_1(\xi) + \beta V_2(\xi) : \xi \in B, \alpha, \beta\in \setR \},
\end{eqnarray*}
where $V_1, V_2\colon B \to TTTM$ are smooth maps
\begin{eqnarray*}
  V_1 = DS \circ \iota, \quad
  V_2 = DS\circ E_1\circ \pi_1 \circ \iota. 
\end{eqnarray*}
Now $\pi_2\circ V_i= \iota$ for $i=1,2$.  A local calculation shows
that $V_1$ and $V_2$ are pointwise linearly independent. Hence
$\Delta$ is a semi-distribution on $TTM$, and by Lemma
\ref{lemma:subman}, $\kappa_3(\Delta)$ is diffeomorphic to $B \times
\setR^2$.

To prove that $\Delta$ is geodesically invariant, let $\gamma\colon
I\to TTM$ be a geodesic of $S^{cc}$ with $\gamma'(0)\in \Delta$. 
By Proposition \ref{prop:NewFromOld}, it follows that
\begin{eqnarray}
\label{eq:parJas}
   \gamma(t) = \kappa_2\circ J'(t), \quad t\in \varepsint
\end{eqnarray}
for a parallel Jacobi field $J\colon \varepsint \to TM$.
By Lemma \ref{lemma:extensionResult} \emph{(ii)} we can extend $J$
into a parallel Jacobi field $J\colon I\to TM$ such that
\eqref{eq:parJas} holds for all $t\in I$. If $t_0\in I$, we have
$\gamma'(t_0) = (\kappa_2\circ \tilde J')'(0)$ for parallel Jacobi
field $\tilde J(t) = J(t+t_0)$, and \emph{(ii)}
follows. Property \emph{(iii)} follows since both submanifolds are
diffeomorphic to $B\times \setR^2$.
\end{proof}

Let us note that configuration space $\pi_2(\Delta)$ 
is a proper subset of $TTM$, and 
$$
  \dim \pi_2(\Delta) = 2n+2, \,\,\,
  \dim \Delta = 2n+2, \,\,\,
  \dim (D\pi_1)(\Delta) = 2n.
$$
The next proposition shows that geodesics in $P$  are in one-to-one
correspondence with parallel Jacobi fields for $M$. 

\newcommand{\confDe}[0]{\operatorname{Conf}(\semiDistribution) }


\begin{proposition}[Geodesics in $\spaceN$] 
\label{prop:projection:Nast}
Suppose $\gamma \colon I\to TTM$ is a curve. Then the following are
equivalent:
\begin{enumerate}
\item $\gamma$ is a geodesic in $\spaceN$.
\item There is a parallel Jacobi field $J\colon I\to TM$ such that
\begin{eqnarray*}
\gamma(t)&=&\kappa_2 \circ J'(t), \quad t\in I.
\end{eqnarray*}
\item There is a geodesic $c\colon I\to M$ and $\alpha, \beta\in
  \setR$ such that
\begin{eqnarray*}
  \gamma(t) = (\alpha + \beta t) c''(t) +
  \beta E_1\circ c'(t), \quad t\in I.
\end{eqnarray*}
\end{enumerate}
Moreover, in (ii) and (iii) $J$ and $c, \alpha, \beta$
are uniquely determined by $\gamma$. 
\end{proposition}

%


The next proposition shows that the geometry of $P$ can be used to
study dynamical properties of $M$.

\begin{proposition} Projection $\pi_{TTM\to M}\colon TTM\to M$ is a
  submersion that maps geodesics in $\spaceN$ into geodesics on $M$.
\end{proposition}



A sub-spray $(S, T^rM, \semiDistribution)$ where $r\ge 0$ is
\emph{complete} if any geodesic $\gamma\colon I\to T^rM$ can be
extended into a geodesic $\gamma\colon \setR\to T^rM$.

\begin{proposition} 
Sub-spray $\spaceN$ is complete if and only if $M$ is complete.
\end{proposition}

\begin{proof}   
  Suppose $P$ is complete. By Proposition \ref{prop:projection:Nast},
  any geodesic $c\colon I\to M$ can be lifted into a geodesic
  $c''\colon I \to TTM$ for $P$. The converse direction follows by
  Proposition \ref{prop:projection:Nast} and Lemma
  \ref{lemma:extensionResult} \emph{(ii)}.
\end{proof}

In general, a geodesic $c\colon I\to M$ for a spray $S$ on $M$ is
uniquely determined by $c'(0)$. The next proposition shows that in
sub-spray $\spaceN$, a geodesic $\gamma \colon I\to TTM$ is uniquely
determined by $\gamma(0)$. This is not surprising in view of 
Proposition \ref{prop:step3a} \emph{(iii)}.

\begin{proposition} 
\label{prop:geoUni1Nast}
  If $\gamma_1\colon I_1 \to TTM$ and $\gamma_2\colon I_2 \to
  TTM$ are geodesics in $\spaceN$ with $\gamma_1(0) = \gamma_2(0)$, then
  $\gamma_1=\gamma_2$ on their common domain.
\end{proposition}

\begin{proof} By Proposition \ref{prop:projection:Nast} we have that
  $\gamma_i = \kappa_2 \circ J_i'$ for parallel Jacobi fields
  $J_i\colon I\to TM$, $i=1,2$. Hence $J_1'(0) = J_2'(0)$, and the
  claim follows.
\end{proof}

Proposition \ref{prop:geoUni1Nast} imposes a strong restriction on the
behavior of geodesics in $P$. For example, if two points in $TTM$ can
be connected with a geodesic in $\spaceN$, then the geodesic is unique
(up to loops). Also, any piece-wise geodesic curve that is continuous
must be smooth.  Therefore $P$ has no broken geodesics nor geodesic
triangles.

For a sub-spray we define conjugate points as for sprays (see
Section \ref{sec:ConjugatePointsForSc}).

\begin{proposition} 
\label{th:mainResult:Nast}
Sub-spray $\spaceN$ has no conjugate points.
\end{proposition}

\begin{proof}
  If a Jacobi field vanishes once, Proposition \ref{prop:geoUni1Nast}
  implies that the corresponding geodesic variation is trivial.
\end{proof}

\subsection*{Acknowledgements}
I.B. has been supported by 
   grant ID 398 from the Romanian Ministry of Education.
M.D. has been supported by 
   Academy of Finland Centre of Excellence programme 213476,
   the Institute of Mathematics at the Helsinki University of Technology,
and
   Tekes project MASIT03 --- Inverse Problems and Reliability of Models. 

%


\providecommand{\bysame}{\leavevmode\hbox to3em{\hrulefill}\thinspace}
\providecommand{\MR}{\relax\ifhmode\unskip\space\fi MR }
\providecommand{\MRhref}[2]{%
  \href{http://www.ams.org/mathscinet-getitem?mr=#1}{#2}
}
\providecommand{\href}[2]{#2}

\end{document}